\magnification=1200
\parindent 0pt
\bigskip

\centerline {\bf Sur les modules d'Iwasawa des
${\bf Z}_p^*/\{\pm 1\}$-extensions.}

\bigskip

\centerline {Marc Perret et Nicolas Saby}

\bigskip

Classification AMS 2000 : primaire 11R23, secondaire 14G15.

\bigskip

Soit $X_{\infty} \cdots \rightarrow X_{n+1} \rightarrow X_n \rightarrow
\cdots X_1 \rightarrow X_0$
une {\it tour} au sens de Mazur-Wiles (voir [12]), c'est \`a dire
une ${\bf Z}_p^*/\{\pm 1\}$-extension d'une courbe alg\'ebrique projective
et lisse
$X_0$ d\'efinie sur le corps fini ${\bf F}_q$ \`a $q$ \'el\'ements de
caract\'eristique $p$,
{\it o\`u chaque point ramifi\'e de
$X_0$ l'est totalement}. Cela signifie que
$X_n$ est un rev\^etement galoisien totalement ramifi\'e de $X_0$, de
groupe de Galois isomorphe
\`a $\left({\bf Z}/p^n{\bf Z}\right)^*/\{\pm1\}$. La limite
$X_\infty = \displaystyle\lim_{\leftarrow} X_n$ est alors un rev\^etement
(ramifi\'e)
infini de $X_0$, de groupe de Galois isomorphe \`a ${\bf Z}_p^*/\{\pm
1\}$.
Nous nous proposons de
donner quelques estimations sur les composantes isotypiques
du {\it module d'Iwasawa} de $X_{\infty} \rightarrow X_0$,
d\'efini de la fa\c con suivante.
Notons $A_n = J_n({\bf F}_q)[p]$ le $p$-Sylow du groupe des points
rationnels
sur ${\bf F}_q$ de la
jacobienne $J_n$ de $X_n$, et $A = \displaystyle \lim_{\leftarrow} A_n$
le {\it module d'Iwasawa} de la ${\bf Z}_p^*/\{\pm1\}$-extension.
Puisque le rev\^etement total est galoisien, le groupe quotient
$Gal(X_1/X_0) \simeq {\bf F}_p^*/\{\pm1\}$ agit sur $A$. On peut donc
d\'ecomposer $A$
en somme directe de {\it composantes isotypiques}

$$A = \sum_{\chi} A^\chi,$$

la somme portant sur les caract\`eres de ${\bf F}_p^*/\{\pm1\}$, c'est
\`a-dire
sur les caract\`eres pairs (i.e. v\'erifiant $\chi(-1) = 1$) de ${\bf
F}_p^*$.
Puisque le groupe du rev\^etement total est ab\'elien, le groupe de Galois
interm\'ediaire $Gal(X_\infty/X_1) \simeq {\bf Z}_p$
agit sur chaque composante $A^\chi$, qui se trouve donc \^etre un
$\Lambda$-module,
o\`u $\Lambda$ d\'esigne {\it l'anneau d'Iwasawa}
$\displaystyle \lim_{\leftarrow} {\bf Z}_p[{\bf Z}/ p^n{\bf Z}] \simeq
{\bf Z}_p[[T]]$.
Le corollaire de la proposition 3 de [12] p. 513 affirme que le
dual de Pontrjagin $M(\chi)$ de
$A^\chi$ est  de type fini et de torsion comme $\Lambda$ module. On peut
donc lui appliquer le
th\'eor\`eme de structure pour les $\Lambda$-modules de type-fini (voir
par exemple [14]) :

\medskip

{\bf Th\'eor\`eme (Iwasawa).} {\it
Si $M$ est un $\Lambda$-module de type-fini, alors il existe un
homomorphisme de $\Lambda$-modules

$$M \longrightarrow
\bigoplus_{i=1}^r \Lambda/{(p^{m_i})} \times \bigoplus_{j=1}^s
\Lambda/{(f_j(T)^{n_j})} \times \Lambda^t$$

de noyau et conoyau finis, o\`u $r, s, t, m_1, \dots, m_r, n_1, \dots ,
n_s$ sont des
entiers positifs, et $f_1, \dots, f_s$ sont des polyn\^omes non constants
d'un certain type de $\Lambda$.

\quad De plus, ces entiers et ces polyn\^omes sont uniquement
d\'etermin\'es par $M$.
}

\bigskip

\quad Les groupes $A_n^\chi= J_n({\bf F}_q)[p]^\chi$ peuvent \^etre
calcul\'es \`a partir de $A^\chi$ :

\medskip

{\bf Proposition (Iwasawa)} {\it
Si les places ramifi\'ees de $X_\infty$ sur $X_0$ le sont totalement,
alors

$$A_n^\chi \simeq A^\chi/((1+T)^{p^n}-1).A^\chi.$$
}

\bigskip

\quad Un calcul \'el\'ementaire montre alors que la croissance du $p$-rang
de $A_n^\chi = J_n({\bf F}_q)[p]^\chi$ est de la forme :

$$d_p(J_n({\bf F}_q)[p]^\chi) =
r_\chi p^{n-1} + c_\chi$$

pour une constante $c_\chi \geq 0$, o\`u $d_p(G)$ d\'esigne le $p$-rang
d'un groupe fini
$G$ et $r_\chi$ est le nombre de facteurs de la forme $\Lambda/p^m$ dans
la d\'ecomposition
d'Iwasawa du dual de Pontrjagin $M(\chi)$.
Une des questions importantes en th\'eorie d'Iwasawa
est l'\'etude des {\it invariants d'Iwasawa} $\mu_\chi$ et $\lambda_\chi$,
o\`u
$\mu_\chi = \sum_{i=1}^{r_\chi} m_i \geq r_\chi$ et
$\lambda_\chi = \sum_{j=1}^{s_\chi} n_j \deg f_j$. En particulier, la
nullit\'e de
$\mu_\chi$ (o\`u, ce qui revient au m\^eme, de $r_\chi$), est d'un grand
inter\^et.
L'objet essentiel de cet article est de prouver le th\'eor\`eme suivant :

\medskip

{\bf Th\'eor\`eme 1.} {\it Soit $X_\infty \rightarrow X_0$
une ${\bf Z}_p^*/\{\pm 1\}$-extension d'une courbe
alg\'ebrique projective irr\'eductible lisse $X_0$ de genre $g(X_0)$
d\'efinie sur le corps fini
${\bf F}_q$. On note $\rho$ (respectivement $\bar\rho$) le nombre de
points ferm\'es (resp. g\'eom\'etriques)
de $X_0$, ramifi\'es dans $X_\infty$.
On suppose que $\rho > 0$, et que chacun de ces points est totalement
ramifi\'e dans $X_\infty$.
Soit $\chi$ un caract\`ere pair de ${\bf F}_p^*$. Alors :

\medskip

(i)  L'invariant $r_\chi$ est major\'e par
$$2g(X_0) + (\rho + \bar\rho) - 1.$$

(ii) L'invariant global
$\lambda = \sum_{\chi} \lambda_\chi$ est minor\'e par $\rho$.}

\bigskip

\quad On peut aussi obtenir une majoration des $r_\chi$ \`a l'aide de la
cohomologie \'etale $p$-adique
en modifiant la technique de Crew dans [3]. Nous en esquisserons les
\'etapes dans la remarque 2 ci-dessous.
L'inter\^et de la preuve que nous proposons ici est quadruple. Elle donne
une majoration meilleure,
elle est parfaitement \'el\'ementaire,
elle donne des r\'esultats pour les corps de nombres (voir la remarque 3),
et enfin elle donne un bien meilleur r\'esultat pour
le caract\`ere trivial et les tours d'Igusa (qui sont les prototypes des
{\it tours}) comme l'affirme le th\'eor\`eme 2
suivant. Commen\c cons par rappeler ce que sont les tours d'Igusa.

\medskip

\quad
Lorsque $N$ est un entier premier \`a $p$,
la courbe d'Igusa $Ig(Np^n)$ de niveau $Np^n$ est
le mod\`ele projectif lisse sur le corps premier ${\bf F}_p$
\`a $p$-\'el\'ements de la courbe affine repr\'esentant (relativement
si $N=1$) le probl\`eme de module
(voir [11], [12] ou [13])

$$S(\hbox{\rm sch\'ema de caract\'eristique} p)
\rightarrow
\left\{\matrix{\hbox{\rm Classes d'isomorphismes} \cr
\hbox{\rm de courbes elliptiques } E/S, \cr
\hbox{\rm avecÊun point }
P \hbox{\rm  d'ordre } N \cr
\hbox{\rm et un g\'en\'erateur } Q  \hbox{\rm de}\cr
\hbox{\rm l'it\'er\'e du Verschiebung } V^n. \cr}\right.$$

\quad  Il s'agit d'un rev\^etement galoisien de $X_1(N)$ de groupe
$\left( {\bf Z}/p^n{\bf Z}\right)^*/\{\pm1\}$ si $N \leq 4$, et
$\left( {\bf Z}/p^n{\bf Z}\right)^*$
sinon, totalement ramifi\'e en les points
supersinguliers de $X_1(N)$, et non ramifi\'es ailleurs. La
limite projective est donc une tour, appell\'ee tour d'Igusa
au dessus de la courbe modulaire $X_1(N)$.

\quad Soit $M(p, \chi)$ la limite projective des duaux de Pontrjagin des
modules
$J_{Ig(p^n)}({\bf F}_p)[p]^\chi$ pour un nombre premier $p$ et un
caract\`ere pair $\chi$
de ${\bf F}_p^*$ donn\'es. Mazur et Wiles appelent ([12], p. 515) {\it
Igusa-r\'egulier}
un tel couple $(p, \chi)$, pour lequel $M(p, \chi) = 0$. Ils prouvent
alors qu'un couple $(p, \chi)$
Igusa-r\'egulier est r\'egulier au sens ordinaire. R\'eciproquement,
ils posent la question ([12], p. 518) de savoir s'il arrive souvent qu'un
couple $(p, \chi)$
soit Igusa-irr\'egulier, bien que classiquement r\'egulier. Le r\'esultat
suivant, qui est une
cons\'equence directe du th\'eor\`eme 1, $(\imath \imath)$ pour $N=1$,
apporte une r\'eponse tr\`es
partielle \`a cette question. Il affirme par exemple que m\^eme si $p$ est
r\'egulier au sens ordinaire,
il existe au moins un caract\`ere $\chi$ pour lequel $(p, \chi)$ soit
Igusa-irr\'egulier.

\medskip

{\bf Corollaire.} {\it Pour tout nombre premier $p$ impair,
il existe au moins un caract\`ere pair $\chi$ de ${\bf F}_p^*$, pour
lequel le couple
$(p, \chi)$ soit Igusa irr\'egulier.}

\bigskip

\quad D'autre part, on peut donner une majoration de la composante
triviale du rang
$r_{N, 1}$ pour ces tours d'Igusa :

\medskip

{\bf Th\'eor\`eme 2.} {\it Pour le caract\`ere trivial $1$ de ${\bf
F}_p^*$, et pour la tour
d'Igusa au dessus de $X_1(N)$, on a la majoration
$$r_{N, 1} \leq s_N$$
o\`u $s_N$ est le nombre de points supersinguliers de $X_1(N)$.}

\bigskip

\quad {\bf Preuve de la premi\`ere assertion du th\'eor\`eme 1.}
D'apr\`es le th\'eor\`eme de structure des $\Lambda$-modules
de type fini d\'ej\`a \'evoqu\'e, pour $n$
assez grand, le $p$-rang $a_n$ de $J_n({\bf F}_q)[p]^{\chi}$ v\'erifie
$a_n = r_{\chi} {p^{n-1}} + c_\chi$. Il existe par la
th\'eorie du corps de classes un rev\^etement non-ramifi\'e $Y_n$ de
$X_n$,
de groupe de Galois isomorphe \`a $\left( {\bf Z}/p{\bf Z} \right)^{a_n}$,
et tel que
$Gal(X_n/X_0)$ agisse sur $Gal(Y_n/X_n)$ par
$\sigma x \sigma^{-1} = x^{\chi(\sigma)}$.
Soit $P$ l'un des points ferm\'es ramifi\'es de
$X_n$. Puisqu'il est non ramifi\'e dans $Y_n$, son
groupe de d\'ecomposition est un sous-groupe cyclique de
$Gal(Y_n/X_n) \simeq \left( {\bf Z}/p{\bf Z} \right)^{a_n}$.
Il est donc trivial ou cyclique d'ordre $p$.
Le rev\^etement ab\'elien maximal interm\'ediaire $Z_n$
non ramifi\'e de $X_n$, de groupe de Galois
de la forme $\left( {\bf Z}/p{\bf Z} \right)^{b_n}$, o\`u $Gal(X_n/X_0)$
agisse selon $\chi$, et o\`u de plus les $\rho$
points ramifi\'es de $X_n$ sont totalement d\'ecompos\'es, v\'erifie donc
$b_n \geq r_\chi {p^{n-1}} + c_\chi - \rho$. Pour fixer les
id\'ees, et par analogie avec le cas des tours d'Igusa,
on appelera {\it supersingulier} les points de $X_n$ (totalement) ramifies
dans la tour.

$$\matrix{
          ~ & Z_n & ~ \cr
      \left. \matrix{\hbox{\rm non ramifi\'e,} \cr
                     \hbox{\rm supersinguliers} \cr
                     \hbox{\rm totalement decompos\'es} \cr
}
                     \right\{
                               &\downarrow &\Biggl) \left({\bf Z}/p{\bf
Z}\right)^{b_n}
     ~~~~~~~~~~\cr
          ~ &X_n &  ~ \cr
     ~~~~~ \left. \matrix{\hbox{\rm supersinguliers} \cr
\hbox{\rm totalement ramifi\'es,}\cr
\hbox{\rm non ramifi\'e ailleurs} \cr
}
\right\{
&\downarrow & \Biggl) \left({\bf Z}/p^n{\bf Z}\right)^*/\{\pm 1\}\cr
&X_0 & ~ \cr
}$$

\quad L'action du groupe cyclique $Gal(X_n/X_0)$
sur $\left({\bf Z}/p{\bf Z}\right)^{b_n}$ par conjugaison est
d\'etermin\'ee par la matrice
$M$ image dans $Aut\left({\bf Z}/p{\bf Z}\right)^{b_n} = GL_{b_n}
\left({\bf Z}/p{\bf Z}\right)$
d'un de ses g\'en\'erateurs. Puisque $M^{{p-1\over 2}p^{n-1}} = I$,
par r\'eduction de Jordan,
$M$ s'\'ecrit dans une base convenable comme une matrice diagonale par
blocs
$M = diag(A_1, \cdots, A_{c_n})$, o\`u chaque $A_i$ est un bloc
\'el\'ementaire de la forme

$$A_i =
\pmatrix{
\zeta & 1 &0&\cdots&0 \cr
0 & \zeta &1&0&0 \cr
\cdots&\cdots&\cdots&\cdots&\cdots&\cr
0&\cdots&0&\zeta&1\cr
0 &\cdots&\cdots& 0&\zeta
\cr},$$

pour une certaine racine ${p-1\over 2}$-i\`eme de l'unit\'e $\zeta$ dans
${\bf F}_p$
commune aux $A_i$, pr\'ecis\'ement puisque l'on ne consid\`ere qu'une
composante isotypique.

\quad Nous affirmons que $A_i$ est de taille $\leq p^{n-1}$. Supposons en
effet le contraire. Le
coefficient sur la premi\`ere ligne et la $(p^{n-1}+1)$-i\`eme colonne de
$I = A_i^{{p-1\over 2}p^{n-1}} = (\zeta I + J)^{{p-1\over 2}p^{n-1}}$
(o\`u $J$ est la matrice avec des z\'eros partout,
sauf sur la seconde diagonale o\`u il y a des 1) serait alors \'egal au
coefficient bin\^omial
$C_{{p-1\over 2}p^{n-1}}^{p^{n-1}}$, qui n'est pas nul modulo $p$, ce qui
est une contradiction !
Par la minoration ci-dessus de $b_n$, on a donc
$c_n.p^{n-1} \geq b_n >r_\chi p^{n-1} + c_\chi - \rho$, c'est \`a dire
(pourvu que $n$ soit assez
grand)

$$c_n \geq r_\chi.$$

\quad Consid\'erons alors la courbe $T_n$ correspondant au quotient
maximal
$\left({\bf Z}/p{\bf Z}\right)^{c_n}$
de \break $Gal(Z_n/X_n)$, sur lequel $Gal(X_n/X_0)$
agisse par {\it homoth\'eties}.
On est en pr\'esence d'une extension

$$1 \rightarrow Gal(T_n/X_n) = \left({\bf Z}/p{\bf Z}\right)^{c_n}
\rightarrow Gal(T_n/X_0)
\rightarrow {\bf Z}/{p-1\over 2}p^{n-1}{\bf Z}\rightarrow 1,$$

o\`u le quotient ${\bf Z}/{p-1\over 2}p^{n-1}{\bf Z}$ agit par $\bar k
\rightarrow \zeta^kId$
sur le noyau ab\'elien
$\left({\bf Z}/p{\bf Z}\right)^{c_n}$.
Soit $\sigma$ un \'el\'ement de $Gal(T_n/X_0)$, dont la restriction
$\bar \sigma$ \`a $X_n$ engendre le quotient
$Gal(X_n/X_0) \simeq {\bf Z}/{p-1\over 2}p^{n-1}{\bf Z}$.

$$
\matrix{
\sigma \in Gal(T_N/X_0)
&
\left(
\matrix{
A \hskip -12mm &
\left(\matrix{
T_n & ~ & ~ \cr
\downarrow & \Bigl) &E \simeq \left( {\bf Z}/p{\bf Z} \right)^{c_n}
~~~~~~~~~~\cr
X_n & ~ & ~ \cr
\downarrow & \Bigl) &<\bar \sigma^{p-1\over 2}> \simeq {\bf Z}/p^{n-1}{\bf
Z} \cr
X_1 & ~ & ~ \cr
}\right.
\cr
~ & ~~~~~~~~~\matrix{
 \downarrow& \Bigl) &<\sigma_{\mid_{X_1}}> \simeq \left({\bf Z}/p{\bf
Z}\right)^*/\{\pm 1\} \cr
 X_0 & ~ & ~ \cr
}
\cr
}\right.
\cr
}$$

\quad Notons $E$ le groupe de Galois $Gal(T_n/X_n)$. Puisque $\sigma$ agit
sur $E$
par $\zeta Id_E$,
$\sigma^{p-1\over 2}$ y agit trivialement. Cela signifie
que le groupe $Gal(T_n/X_1)$ est un $p$-groupe {\it ab\'elien} $A$.
On consid\`ere alors l'un des $\rho$
points ferm\'es ``supersinguliers'' de $X_1$. Puisqu'il est totalement
ramifi\'e dans l'extension
cyclique interm\'ediaire
$X_n$, puis totalement d\'ecompos\'e dans $T_n$, il a un groupe de
d\'ecomposition cyclique d'ordre $p^{n-1}$ dans $Gal(T_n/X_1)$. Les
groupes de d\'ecompositions des points supersinguliers de $X_1$ engendrent
donc un sous-groupe
-que l'on notera $D$- de $A = Gal(T_n/X_1)$, engendr\'e par $\rho$
\'el\'ements.
La courbe $U_n$ fix\'ee par $D$
est donc un rev\^etement $p$-\'el\'ementaire de $X_1$ {\it non ramifi\'e}
(puisque $T_n$ n'est ramifi\'ee sur $X_0$ qu'en les supersinguliers), de
groupe
de Galois isomorphe \`a $\left( {\bf Z}/p{\bf Z} \right)^{d_n}$ avec

$$d_n \geq c_n-\rho \geq r_\chi  - \rho.$$

\quad On en d\'eduit que $U_n$ est un rev\^etement {\it mod\'er\'ement
ramifi\'e} de $X_0$,
ramifi\'e seulement en les supersinguliers, et o\`u au moins un point
ferm\'e est
totalement d\'ecompos\'e. Par le th\'eor\`eme de
Grothendieck (voir [5], chapitre XIII, cor 2.12), le groupe de Galois
$Gal(U_n/X_0)$ peut donc \^etre engendr\'e par $2g(X_0) + \bar\rho - 1$
\'el\'ements.

\quad Mais ce groupe est par le th\'eor\`eme de Burnside le produit
semi-direct de
$Gal(U_n/X_1)$ par $Gal(X_1/X_0)$.
Un \'el\'ement de $Gal(U_n/X_0)$ peut donc \^etre repr\'esent\'e par un
couple
$(x,  \sigma_{\mid_{X_1}}^k) \in \left( {\bf Z}/p{\bf Z}
\right)^{d_n}\times Gal(X_1/X_0)$,
la loi de groupe \'etant

$$(x, \sigma_{\mid_{X_1}}^k)\times (y,  \sigma_{\mid_{X_1}}^\ell)
= (x+\zeta^k y,  \sigma_{\mid_{X_1}}^{k+\ell}).$$

\quad On observe que le projet\'e sur le premier facteur d'un produit est
dans
l'espace vectoriel engendr\'e par les projet\'es des deux termes. En
particulier,
toute partie g\'en\'eratrice de $Gal(U_n/X_0)$ a au moins
$d_n$ \'el\'ements. Il en r\'esulte que

$$r_\chi - \rho \leq d_n \leq 2g(X_0) + \bar\rho - 1,$$

ce qui n'est autre que l'assertion {\it (i)}.

\bigskip

\quad {\bf Preuve du th\'eor\`eme 2.} Supposons que cette majoration soit
fausse, et reprenons mot pour mot
la preuve pr\'ec\'edente jusqu'\`a
l'introduction de la courbe $T_n$. Dans le cas qui nous occupe,
$X_0 = X_1(N)$ est la r\'eduite modulo $p$ de la courbe modulaire pour le
sous-groupe de congruence $\Gamma_1(N)$,
$X_1$ est la courbe d'Igusa $Ig(Np)$ de niveau $Np$ et $\rho$ est le
nombre de points supersinguliers
de $X_1(N)$, traditionnellement not\'e $\rho = s_N$. Puisqu'on consid\`ere
le caract\`ere trivial,
$Gal(T_n/X_0)$
est alors extension d'un groupe ab\'elien par un groupe cyclique
{\it avec action triviale}. Il
est donc lui m\^eme ab\'elien, ce qui implique qu'il ne peut \^etre
engendr\'e par moins
de $c_n$ \'el\'ements (avec toujours $c_n \geq r_{N, 1}$). D'autre part,
les
$s_N$ points supersinguliers de $X_0 = X_1(N)$ ont un groupe de
d\'ecomposition cyclique
dans $Gal(T_n/X_0)$, puisqu'ils sont totalement ramifi\'es dans
l'extension cyclique $X_n/X_0$,
puis totalement d\'ecompos\'es dans $T_n/X_n$. Le sous-groupe $D_1$ de
$Gal(T_n/X_0)$
engendr\'e par les supersinguliers de $X_1(N)$ est donc strictement
contenu dans
$Gal(T_n/X_0)$ puisqu'on a suppos\'e que la majoration \'etait fausse.
La courbe fix\'ee par $D_1$ est donc un rev\^etement
{\it non ramifi\'e} de $X_1(N)$, {\it non-trivial},
o\`u les supersinguliers de $X_1(N)$ sont
{\it totalement d\'ecompos\'es}. Soit $X(N)$ la r\'eduite modulo
$p$ de la courbe modulaire pour le sous-groupe de congruence $\Gamma(N)$.
Par produit fibr\'e avec $X(N)$ au dessus de $X_1(N)$,
cela fournit un rev\^etement de $X(N)$ du m\^eme type, toujours non
trivial puisque $X(N)$
est de degr\'e $N$ premier \`a $p$ sur $X_1(N)$. Cela est en contradiction
avec le
th\'eor\`eme suivant :

\medskip

{\bf Th\'eor\`eme (Ihara, [8]).} {\it Soit $N$ un entier premier \`a $p$.
Alors $X(N)$ ne poss\`ede
pas de rev\^etement galoisien non ramifi\'e et non trivial, o\`u les
points supersinguliers soient
totalement d\'ecompos\'es.}

\bigskip

\quad {\bf Preuve de la seconde assertion du th\'eor\`eme 1.} Soit $n \geq
1$. On note pour simplifier $G_n$
le groupe de Galois de $X_n$ sur $X_0$ et $P(X_n)$ le groupe des diviseurs
principaux sur $X_n$.
La suite exacte

$$1 \longrightarrow P(X_n) \longrightarrow Div^0(X_n) \longrightarrow
J_n({\bf F}_q)
\longrightarrow 1$$

donne

$$1 \longrightarrow P(X_n)^{G_n} \longrightarrow Div^0(X_n)^{G_n}
\longrightarrow J_n({\bf F}_q)^{G_n},$$

d'o\`u une compos\'ee de deux injections

$$Div^0(X_n)^{G_n}/P(X_n)^{G_n} \longrightarrow J_n({\bf F}_q)^{G_n}
\longrightarrow J_n({\bf F}_q).\leqno (1)$$

\quad D'autre-part il y a deux suites exactes

$$1 \longrightarrow P(X_n)^{G_n}/P(X_0) \longrightarrow
Div^0(X_n)^{G_n}/P(X_0)
\longrightarrow Div^0(X_n)^{G_n}/P(X_n)^{G_n}
\longrightarrow 1 \leqno (2)$$
et

$$ Div^0(X_n)^{G_n}/P(X_0)
\longrightarrow Div^0(X_n)^{G_n}/Div^0(X_0)
\longrightarrow 1. \leqno (3)$$

\quad Mais :

\quad 1) Puisque
$Div^0(X_n)^{G_n}/Div^0(X_0) \simeq \prod_{P\in X_0 ~ram~dans~X_n}{\bf
Z}/e_P{\bf Z}$, et puisque
$X_n \longrightarrow X_0$ est totalement ramifi\'e en
$\rho$ points avec
indice de ramification $e = {p-1\over 2}p^{n-1}$,
on obtient, par (3), que

$$\left({\bf Z}/{p-1\over 2}p^{n-1}{\bf Z}\right)^\rho
est~un~quotient~de~Div^0(X_n)^{G_n}/P(X_0).
\leqno (4)$$

\medskip

\quad 2) D'autre part, la suite exacte

$$1 \longrightarrow {\bf F}_q^* \longrightarrow {\bf F}_q(X_n)^*
\longrightarrow P(X_n)
\longrightarrow 1$$

donne

$$1 \longrightarrow {\bf F}_q^* \longrightarrow {\bf F}_q(X_0)^*
\longrightarrow P(X_n)^{G_n}
\longrightarrow H^1(G_n,{\bf F}_q^*) \longrightarrow H^1(G_n,{\bf
F}_q(X_n)^*) = 1,$$

c'est \`a dire

$$1 \longrightarrow P(X_0) \longrightarrow  P(X_n)^{G_n}
\longrightarrow H^1(G_n,{\bf F}_q^*) \longrightarrow1.$$

\quad Le $H^1$ en jeu est calcul\'e pour l'action triviale de $G_n$ sur
${\bf F}_q^*$, donc

$$H^1(G_n,{\bf F}_q^*) = Hom(G_n, {\bf F}_q^*).$$

Puisque
$G_n = {\bf Z}/{p-1\over 2}p^{n-1}{\bf Z}$, on a

$$P(X_n)^{G_n}/P(X_0) = H^1(G_n,{\bf F}_q^*) \simeq End({\bf Z}/{p-1\over
2}{\bf Z})
~\hbox{\rm ne d\'epend pas de }Ên. \leqno (5)$$

\medskip

\quad Lorsque l'on injecte ces deux points (4) et (5) dans la suite exacte
(2),
on obtient que, \`a un groupe
born\'e pr\`es, $Div^0(X_n)^{G_n}/P(X_n)^{G_n}$, et donc aussi $J_n({\bf
F}_q)$,
est au moins aussi gros que
$\left({\bf Z}/{p-1\over 2}p^{n-1}{\bf Z}\right)^\rho$.
Mais il est ais\'e de s'assurer, \`a l'aide de la proposition d'Iwasawa
rappel\'ee dans l'introduction,
que l'invariant global $\lambda$ n'est autre
que le nombre de facteurs cycliques du $p$-Sylow de $J_n({\bf F}_q)$ dont
les ordres tendent vers l'infini avec $n$.
Il en r\'esulte bien que $\lambda$ est minor\'e par $\rho$.

\bigskip

{\bf Remarque 1. Tours de Shimura et tours de Mazur-Wiles.}
Le th\'eor\`eme 2 ne se
g\'en\'eralise pas au cas g\'en\'eral des tours faute d'un
th\'eor\`eme \`a la Ihara.
Cependant,
toujours d'apr\`es Ihara, un tel th\'eor\`eme
existe pour les courbes de Shimura (voir [9]): \'etant donn\'e une courbe
de Shimura $S$ sur ${\bf F}_p$,
il existe un groupe fini $G_S$, tel que tout rev\^etement galoisien X de
$S$ non ramifi\'e,
o\`u les points supersinguliers sont totalement d\'ecompos\'es, ait pour
groupe de Galois
un quotient de $G_S$. Par exemple, le th\'eor\`eme de Ihara \'enonc\'e
lors de la preuve du th\'eor\`eme 2 affirme que pour la courbe modulaire
$X(N)$ de niveau
$N$ premier \`a $p$ r\'eduite modulo $p$, on a
$G_{X(N)} = 1$.
On peut fabriquer une tour au sens de Mazur-Wiles d'Igusa-Shimura
\`a partir des composantes d'Igusa des courbes de Shimura
(Cf [1] ou [2]).
 On pourrait donc donner une majoration de l'invariant $r$
d'une tour d'Igusa-Shimura pour le caract\`ere trivial en
fonction du $p$-rang de $G_S$.
Malheureusement, la preuve de l'existence de $G_S$ est totalement
ineffective, ce qui
\^ote tout inter\^et \`a une telle majoration.

\bigskip

{\bf Remarque 2. Majoration de $r_\chi$ \`a l'aide de la cohomologie
$p$-adique.} On ne reprend
que partiellement les notations de [3] pages 34 et 35 :
soit $S/k$ un sch\'ema s\'epar\'e de type fini sur un corps
alg\'ebriquement clos
$k$ de caract\'eristique $p$ et $G$ un $p$-groupe fini agissant librement
sur $S$, de quotient $S/G$.
On suppose qu'un groupe
$\Delta$ fini d'ordre premier \`a $p$ agit aussi sur $S$. Il agit donc sur
la suite
spectrale de Hochschild-Serre. Si les actions de $G$ et de $\Delta$
commutent, il y a alors pour chaque caract\`ere
$\chi$ de $\Delta$
une suite spectrale pour les composantes isotypique $H^i_c(S, {\bf
Q}_p)^\chi$, qui montre que

$$H^i_c(S/G, {\bf Q}_p)^\chi = (H^i_c(S, {\bf Q}_p)^\chi)^G.$$

\quad Supposons maintenant que $Y$ soit rev\^etement
ramifi\'e d'une courbe alg\'ebrique projective est lisse $X$ sur $k$,
galoisien de groupe de Galois un $p$-groupe $G$, et que $\Delta$ agisse
sur $X$ en commutant \`a $G$.
Si $C$ est l'une de ces deux courbes, la $\chi$-partie du $p$-rang de $C$
est d\'efinie par
$p_C^\chi := dim_{{\bf Q}_p} H^i_c(C, {\bf Q}_p)^\chi$.
La suite exacte d'excision jointe au th\'eor\`eme 1.5 de Crew et \`a la
remarque pr\'ec\'edente montre alors,
de la m\^eme fa\c con que l'on montre le corollaire 1.8 de [3] :

\medskip

{\bf Proposition 1.} {\it Dans la situation pr\'ec\'edente :

\medskip

$(\imath)$ si $\chi$ n'est pas le caract\`ere trivial, alors

$$p_X^\chi = (\sharp G) p_Y^\chi.$$

$(\imath \imath)$ Pour le caract\`ere trivial $1$, on a

$$p_X^1 -1= (\sharp G)(p_Y^1-1) + \sum_{x \in X^{ram}} (e_x - 1).$$
}

\bigskip

\quad On rappelle que
si $C$ est de genre $g_C$, on a
$p_C^\chi \leq p_C \leq g_C$.
Revenons alors \`a notre situation des tours $X_\infty \rightarrow X_0$
sur le corps premier ${\bf F}_p$,
et consid\'erons un caract\`ere non trivial
$\chi$ de $\Delta = Gal(X_1/X_0) = {\bf F}_p^*/\{\pm1\}$.
La proposition 1 appliqu\'ee au rev\^etement ramifi\'e $X_n$ de $X_1$,
puis la
formule de Riemann-Hurwitz appliqu\'ee au rev\^etement {\it mod\'er\'ement
ramifi\'e} $X_1$ de $X_0$
donnent, en supposant pour simplifier que $p$ est impair :

$$r_\chi p^{n-1} + c_\chi \leq p_{X_n}^\chi
= p^{n-1}p_{X_1}^\chi
\leq p^{n-1}g(X_1)
$$

avec $2g(X_1)-2 = {p-1 \over 2}(2g(X_0)-2) + \bar\rho({p-1 \over 2} - 1)$.
On obtient donc la majoration

$$r_\chi \leq{p-1 \over 2}g(X_0) + {p - 3 \over 4}(\bar\rho-2),$$

qui est moins bonne que celle du th\'eor\`eme 1 d\`es que $p>5$.

\bigskip

{\bf Remarque 3. Le cas des corps de nombres.}
La technique de la preuve du th\'eor\`eme 2 fonctionne
tr\`es bien dans le cas des corps de nombres.
En revanche, la preuve de {\it (ii)} ne s'adapte pas du tout \`a ce cas.
Pour ce qui concerne la preuve de {\it (i)}, il manque malheureusement un
analogue
du th\'eor\`eme de Grothendieck.
Cet analogue a cependant \'et\'e conjectur\'e par Harbater
(dans [6] dans le cas o\`u $K = {\bf Q}$, et dans [7] dans le cas
g\'en\'eral) :
pour un corps de nombres $K$ et un id\'eal entier $I$ de $K$, notons
$\pi_A^t(K, I)$ le {\it groupe fondamental
alg\'ebrique mod\'er\'e de $Spec {\cal O}_K[{I^{-1}}]$}, c'est \`a dire
l'ensemble
des groupes finis
apparaissant comme groupe de Galois d'une extension mod\'er\'ement
ramifi\'ee de
${\bf Q}$ et non ramifi\'ee hors de $I$ .

\medskip

{\bf Conjecture. (Harbater)} {\it Il existe une constante absolue $C$,
telle que pour tout corps de nombres $K$ de genre $g(K)$ (au sens de
Weil), ayant
$r(K)$ plongements r\'eels, et pour tout id\'eal entier $I$ de $K$,
tout \'el\'ement de $\pi_A^t(K, I)$ peut-\^etre engendr\'e par
$ 2g(K) + r(K) - 1 + C + \log N_{K/{\bf Q}} (I)$ \'elements.}

\bigskip

\quad Nous rappelons que le groupe de Galois de l'extension
$\displaystyle \bigcup_nK\left( \exp ({2\imath \pi \over p^n}) \right) /K$
est de la forme ${\bf F}_p^*/G \times {\bf Z}_p$ pour un sous-groupe $G$
de ${\bf F}_p^*$
convenable.
\medskip

{\bf Th\'eor\`eme 3.} {\it Soit $p$ un nombre premier impair, $K$ un corps
de nombres,
$G$ un sous-groupe fini de ${\bf F}_p^*$ et
$K = K_0 \subset K_1 \subset \cdots \subset K_{\infty}$ une
${\bf F}_p^*/G \times {\bf Z}_p$-extension de $K$. On suppose que
les $\rho$ id\'eaux
premiers ${\cal P}_1, \cdots, {\cal P}_\rho$ de $K$ au dessus de $p$ sont
totalement ramifi\'es dans
$K_{\infty}$. Soit $\chi$ un caract\`ere de ${\bf F}_p^*/G$. Alors

1. Sous la conjecture d'Harbater, les invariants $r_{p, \chi}$ sont
major\'es par

$$2g(K) + (\rho + [K : {\bf Q}]\log p) -1 + (r(K) + C).$$

2. Pour le caract\`ere trivial $\chi = 1$, l'invariant $r_{p, 1}$ est
major\'e par

$$\rho + d_p \left( {\cal C}l(K)/<\overline {{\cal P}_1}, \cdots,
\overline {{\cal P}_\rho}> \right).$$
}

\bigskip

\quad Comparons ce th\'eor\`eme 3 aux th\'eor\`emes 1 et 2. Dans
l'assertion 1 du th\'eor\`eme 3,
la quantit\'e $[K : {\bf Q}]\log p = \log N_{[K : {\bf Q}]}(p{\cal O}_K)$
est l'analogue pour les
corps de nombres du degr\'e du diviseur r\'eduit de ramification
$\bar\rho = \deg (\sum_{i=1}^\rho P_i) = \sum_{i=1}^\rho\deg  P_i$ de la
situation g\'eom\'etrique.
Dans cette m\^eme assertion, le terme suppl\'ementaire $r(K) + C$
provient de la conjecture d'Harbater. Quant \`a l'assertion 2, le terme
correctif
$d_p \left( {\cal C}l(K)/<\overline {{\cal P}_1}, \cdots,
\overline {{\cal P}_\rho}> \right)$ est l'analogue du $p$-rang du groupe
$G_S$ dont il a \'et\'e
question dans la remarque 1 pour les courbes de Shimura.

\quad Notons que pour $p = 2$, un \'enonc\'e analogue vaut en adjoignant
la racine
quatri\`eme de l'unit\'e $\imath$ \`a $K$.
Par exemple, la seconde assertion montre que les invariants $r_p$ de la
${\bf Z}_p$-extension cyclotomique de
${\bf Q}[\root 3 \of 2]$ sont
inf\'erieurs \`a $3$ pour tout $p$.
Rappelons que, d'apr\`es Ferrero et Washington (voir [4]), les invariants
$\mu_p$
(o\`u, ce qui revient au m\^eme, les invariants $r_p$)
des ${\bf Z}_p$-extensions cyclotomiques des extensions ab\'eliennes de
${\bf Q}$ sont
nuls. Cette nullit\'e a d'ailleurs \'et\'e conjectur\'ee par Iwasawa pour
tous les
corps de nombres. D'autre part, il existe d'apr\`es Iwasawa
(voir [10]) des ${\bf Z}_p$-extensions (non cyclotomiques !) de corps de
nombres, dont
les invariants $\mu$ (et donc aussi $r$) sont non nuls.

\bigskip

{\bf Remerciements.} Nous tenons \`a remercier le rapporteur pour ses
remarques tr\`es utiles.

~

\bigskip

\centerline {\bf Bibliographie.}

\medskip

[1] K. Buzzard, Integral models of certain Shimura curves, {\it Duke Math.
J.} {\bf 87} (1997), no. 3, 591--612.

\medskip

[2] H. Carayol, Sur la mauvaise r\'eduction des courbes de Shimura,
{\it Compositio Math.} {\bf 59} (1986), no. 2, 151--230.

\medskip

[3] R. Crew, Etale $p$-covers in characteristic $p$, {\it Compositio
Math.} {\bf 52} (1984),
31-45.

\medskip

[4] B. Ferrero, L. Washington, The Iwasawa invariants $\mu_p$ vanishes for
abelian number fields, {\it Annals of Maths} {\bf 109} (1979), 377-395.

\medskip

[5] A. Grothendieck, {\it Rev\^etement \'etales et groupe fondamental} SGA
I,
Lect. Notes in Maths. {\bf 224}, Springer, 1971.

\medskip

[6] D. Harbater, "Galois groups with prescribed ramification", in
Arithmetic Geometry,
Childress et Jones editors, Contemporary Mathematics 174 (1994), p. 35-60.

\medskip

[7] D. Harbater, communication priv\'ee.
\medskip

[8] Y. Ihara, "On modular curves over finite fields", in
discrete subgroups of Lie groups and applications to moduli, Oxford 1975,
161-202.

\medskip
[9] Y. Ihara, Congruence relations and fundamental groups, {\it J. of
Algebra} {\bf 75}
(1982), 445-451.

\medskip

[10] K. Iwasawa, "On the $\mu$-invariants of ${\bf Z}_\ell$-extensions",
in
{\it Number Theory, Algebraic Geometry and Commutative algebra} in honnor
to
Y. Akizuki, Konokuniya, Tokyo (1973), 1-11.

\medskip

[11] N. Katz et B. Mazur, {\it Arithmetic moduli of elliptic curves}
Annals of Math. studies,
vol. 108, Princeton Univ. Press, 1985.

\medskip

[12] B. Mazur et A. Wiles, Analogies between function fields and number
fields, {\it Amer. J. Math.}
{\bf 105} (1983), 507-521.

\medskip

[13] N. Saby, Th\'eorie d'Iwasawa g\'eom\'etrique : un th\'eor\`eme de
comparaison,
{\it J. Numb. Th.} {\bf 59} no 2 (1996), 225-247.

\medskip

[14] L. Washington, {\it Introduction to cyclotomic fields}, Graduate
Texts in Mathematics,
Springer-Verlag, New-York, 1982.

\bigskip

Marc PERRET

GRIMM

Universit\'e de Toulouse II le Mirail

5 All\'ees Antonio Machado

31 058 TOULOUSE Cedex

FRANCE

perret@univ-tlse2.fr

\bigskip

~
\bigskip

Nicolas SABY

D\'epartement de Math\'ematiques

Universit\'e de Montpellier II

case 051

Place eug\`ene Bataillon

34 095 Montpellier cedex 5

FRANCE

saby@math.univ-montp2.fr

\end